\documentclass[a4paper]{article}
\usepackage{amsthm,amssymb,amsmath,enumerate,graphicx,psfrag,color}

\newtheorem{definition}{Definition}

\newtheorem{theorem}[definition]{Theorem}

\newtheorem{lemma}[definition]{Lemma}
\newtheorem{problem}[definition]{Problem}
\newtheorem{example}[definition]{Example}

\newcommand{\minG}{\delta^{V,\Omega}(G)}
\newcommand{\comment}[1]{}
\newcommand{\N}{\mathbb N}

\newcommand{\df}{d_{e/v}}

\title{Forcing large complete minors in infinite graphs}
\author{Maya Stein\footnote{Centro de Modelamiento Matem\'atico, Universidad de Chile, Santiago, Chile, email: {\tt mstein@dim.uchile.cl}, Fondecyt grant no. 1104099} \ \  \& Jos\'e Zamora\footnote{Universidad Andr\'es Bello, Santiago, Chile}}
\date{30.12.2010}
\usepackage{graphicx}

\begin{document}
\maketitle

\begin{abstract}
It is well-known that in finite graphs, large complete minors/topological minors can be forced by assuming a large average degree. 
Our aim is to extend this fact to infinite graphs. For this, we generalise the notion of the relative end degree, which had been previously introduced by the first author for locally finite graphs, and show that large minimum relative degree at the ends and large minimum degree at the vertices imply the existence of large complete (topological) minors in infinite graphs with countably many ends. \end{abstract}

\section{Introduction}
A recurrent question in finite graph theory is how certain substructures, such as specific minors or subgraphs, can be forced by density assumptions. The density assumptions are often expressed via a lower bound on the average or minimum degree of the graph. 

The classical example for this type of questions is Tur\'an's theorem, or, in the same vein, the Erd\H os-Stone theorem. Infinite analogues of these results are not difficult: it is easy to see that if the upper density of an infinite graph $G$ is at least the Tur\'an density for $k$, i.e.~$(k-2)/(k-1)$, then $K^ k\subseteq G$, that is, $G$ contains  the complete graph of order $k$ as a subgraph (see Bollob\'as~\cite{BBModern}, and also~\cite{mayaMCW2010}). 
On the other extreme, a result of Mader affirms that an average degree of at least $4k$ ensures the existence of a $(k+1)$-connected subgraph. This result has been extended to infinite graphs by the first author~\cite{hcs}. 

Half-way between the two types of results just discussed lie two well-known results due to Kostochka~\cite{kostochka1984} and to Bollob\'as and Thomason~\cite{bollobas1998}, respectively, which we sum up in the following  theorem. 

\begin{theorem}\label{thm:minor}$\!\!${\bf\cite{DBook}}
 There are constants $c_1,c_2$ so that for each  $k\in\N$, and each graph $G$ the following holds. If $G$ has average degree at least $c_1 k\sqrt{\log k}$ then $K^k\preceq G$ and if $G$ has average degree at least $c_2 k^2$ then $K^k\preceq_{top} G$.
\end{theorem}

Our aim is to extend this result to infinite graphs (qualitatively, that is, without necessarily using the same functions $c_1 k\sqrt{\log k}$ and $c_2 k^2$ from Theorem~\ref{thm:minor}). We call a finite graph $H$ a minor/topological minor of a graph $G$ if  it is a minor/topological minor of some finite subgraph of $G$.\footnote{Note that this means all branch-sets of our minor are finite, but as long as $H$ is finite, it clearly  makes no difference whether we allow infinite branch-sets or not.}

Avoiding the difficulty of defining an average degree for infinite graphs (the upper density mentioned above is too strong for our purposes, cf.~\cite{mayaMCW2010}), we shall stick to the minimum degree for our extension of Theorem~\ref{thm:minor} to infinite graphs. This works fine for rayless graphs: the first author showed in~\cite{InfExt} that every rayless graph of minimum degree $\geq m\in\N$ has a finite subgraph of minimum degree $\geq m$. Thus, Theorem~\ref{thm:minor} can be literally extended to rayless graphs if we replace the average with the minimum degree.

In general, however, large minimum degree at the vertices alone is not strong enough to force large complete minors. This is so because of the infinite trees, which may attain any minimum degree condition without containing any interesting substructure. So we need some additional condition that prevents the density from `escaping to infinity'.

The most natural way to impose such an additional condition is to impose it on the ends\footnote{The {\em ends} of a graph are the equivalence classes of the rays (i.e.~$1$-way infinite paths), under the following equivalence relation. Two rays are equivalent if no finite set of vertices separates them. (A set $S$ of vertices is said to {\em separate} two rays $R_1$, $R_2$ if $V(R_1)\setminus S$, and $V(R_2)\setminus S$ lie in different components of $G-S$.) The set of ends of a graph $G$ is denoted by $\Omega(G)$. See the infinite chapter of~\cite{DBook} for more on ends.} of the graph. This  approach has also proved successful in other
recent work~\cite{degree, InfExt,hcs}.
In this way, that is, defining the degree of an end in appropriate way, the minimum degree, now taken over vertices and ends, can continue to serve as our condition for forcing large complete minors.

In~\cite{InfExt}, the  {\em relative degree} of an end was introduced for locally finite graphs. 
In order to explain it, a few notions come in handy. 

Let $G$ be a locally finite graph.
 The {\em edge-boundary} of a subgraph $H$ of $G$, denoted by $\partial^G_e H$, or $\partial_e H$ where no confusion is likely, is the set $E(H,G-H)$.
The {\em vertex-boundary} $\partial^G_v H$, or $\partial_v H$, of a subgraph $H$ of $G$ is the set of all vertices in $H$ that have neighbours in $G-H$. 
Now, the idea of the relative degree is to calculate  the ratios of $|\partial_eH_i|/|\partial_vH_i|$  of certain subgraphs $H_i$ of $G$,  and then take the relative degree to be the limit of these ratios as the $H_i$ in some sense converge to $\omega$. 

For this, define an $\omega$-region of an end $\omega$ of a graph $G$ as an induced connected subgraph which contains some ray of $\omega$ and whose vertex-boundary is finite. 
For $V'\subseteq V(G)$, $\Omega'\subseteq \Omega(G)$, a $V'$--$\Omega'$ separator is a set $S\subseteq V(G)$ such that $V'\not\subseteq S$ and  such that no component of $G\setminus S$ contains both a vertex of $V'$ and a ray from $\omega$.  We write  $\Omega^G(H)$ for the set of all ends of $G$ that have a ray in $H\subseteq G$.

  Write $(H_i)_{i\in\N}\rightarrow\omega$ if $(H_i)_{i\in\N}$ is an infinite sequence of distinct $\omega$-regions  such that $H_{i+1}\subseteq H_i- \partial_v H_i$ and $\partial_v H_{i+1}$ is a minimal $\partial_v H_i$--$\Omega^G (H_{i+1})$ separator in $G$, for each $i\in\N$. Now, define 
\[
\df(\omega):=\inf_{(H_i)_{i\in\N}\rightarrow\omega}\liminf_{i\rightarrow\infty}\frac{|\partial_e H_i|}{|\partial_v H_i|}.
\]
Note that it does not matter whether we consider the $\liminf$ or the $\limsup$, because if $(H_i)_{i\in\N}\rightarrow\omega$, also  all subsequences of $(H_i)_{i\in\N}$ converge to $\omega$. (For the same reason we could restrict our attention only to sequences $(H_i)$ for which $\lim_{i\rightarrow\infty}\frac{|\partial_e H_i|}{|\partial_v H_i|}$ exists.)

For more discussion of this notion see Section~\ref{sec:disc}. From now on we write $\minG$ for the minimum degree/relative degree, taken over all vertices and ends of  $G$.
The first author showed:
\begin{theorem}$\!\!${\bf\cite{InfExt}}\label{locfin}
Let $m\in\mathbb N$ and let $G$ be a locally finite graph. If $\minG \geq m$, then $G$ has a finite subgraph of average degree at least $m$.
\end{theorem}
In particular, this means that if  $\minG \geq c_1k\sqrt{\log k}$ for a locally finite $G$, then $K^k\preceq G$, and if $\minG \geq c_2 k^2$, then $K^k\preceq_{top} G$. (The $c_i$ are the constants from Theorem~\ref{thm:minor}.)

\medskip

In this paper, we extend the notion  of the relative degree  to arbitrary infinite graphs. For a given end $\omega$ of some infinite graph $G$, let $Dom (\omega)$ denote the set of all vertices that dominate\footnote{We say a vertex $v$ {\em dominates} an end $\omega$ if $v$ dominates some ray $R$ of $\omega$, that is, if there are infinitely many $v$--$V(R)$-paths that are disjoint except in $v$. (It is not difficult to see that then $v$ dominates all rays of~$\omega$.)} $\omega$. If  $G_\omega:=G-Dom(\omega)$ does not contain any rays from $\omega$, or if  $|Dom(\omega)|\geq\aleph_0$, then we set $\df(\omega):=|Dom(\omega)|$. 
Otherwise,  writing $\hat\omega$ for the unique\footnote{The uniqueness of $\hat \omega$ follows at once from the fact that $|Dom(\omega)|<\aleph_0$.} end of $G_\omega$ that contains rays from $\omega$ we define 
\[
\df(\omega):=|Dom(\omega)|+\inf_{(\hat H_i)_{i\in\N}\rightarrow\hat\omega}\liminf_{i\rightarrow\infty}\frac{|\partial_e \hat H_i|}{|\partial_v \hat H_i|}.
\]
Note that  the $\hat H_i$ are $\hat\omega$-regions of $G_\omega$. Also note that for locally finite graphs, our definition coincides with the one given earlier.
For further discussion of our notion of the relative degree, for examples, and for alternative definitions that do (not) work, see Section~\ref{sec:disc}.

Our main result is the following version of Theorem~\ref{locfin} for graphs with countably many ends.

\begin{theorem}\label{thm:cntble}
Let $k\in \N$, let $m\in\mathbb Q$, and let $G$ be a graph such that $|\Omega(G)|\leq\aleph_0$ and $\minG> m$. Then  $K^{k}\preceq_{top}G$ or $G$ has a finite subgraph of average degree greater than $m-k+1$.
\end{theorem}
 
 This result, together with Theorem~\ref{thm:minor}, at once implies the desired extension of Theorem~\ref{thm:minor} to graphs with countably many ends.
 
\begin{theorem}\label{thm:cntbleMinor}
Let $k\in \N$, and let $G$ be a graph with  $|\Omega(G)|\leq\aleph_0$.
\begin{enumerate}[(a)]
\item If $\minG \geq c_1k\sqrt{\log k}+k$, then $K^k\preceq G$. 
\item If $\minG \geq c_2 k^2+k$, then $K^k\preceq_{top} G$. 
\end{enumerate}
\end{theorem}

Theorem~\ref{thm:cntbleMinor} is best possible in the following sense. There is  no function $f:\N\to\N$ so that  $K^k\preceq G$ for every graph with  $\minG \geq f(k)$. To see this, consider the graph we obtain by the following procedure. Let $G_1$ be a double-ray and let $E_1:=E(G_1)$. Now, for $i\geq 2$, replace each edge $e\in E_{i-1}$ with $\aleph_0$ many paths of length two whose middle points are connected by a (fresh) double-ray $D_e$. Let $E_i$ be the edges on the $D_e$. After $\aleph_0$ many steps, we obtain a planar graph $G$ with $\minG =\aleph_0$. Beause of its planarity, $G$ has no  complete minors of order greater than~$4$.

\begin{problem}\label{problem}
Is there a `degree condition'  which forces large complete (topological) minors in arbitrary infinite graphs?
\end{problem}
%

\section{Discussion of the relative degree}\label{sec:disc}

In this section, we will discuss our definition as well as possible alternative definitions of the relative degree of an end. This motivation is not necessary for the understanding of the rest of the paper and may be skipped at a first reading.

\subsection{Large vertex-degree is not enough}

As we have already seen in the introduction, large minimum degree at the vertices alone is not sufficient for forcing large complete (topological) minors in infinite graphs, because of the trees. A similar example discards the following alternative. The  {\em vertex-/edge-degree} of an end $\omega$ was introduced in~\cite{hcs}, see also~\cite{DBook}, as the supremum of the cardinalities of sets of vertex-/edge-disjoint rays in $\omega$. Clearly, the vertex-degree of an end is always at most its edge-degree, so for our purposes we may restrict our attention to the vertex-degree. It is not difficult to show~\cite{arbo} that an end has vertex-degree $\geq k$ if and only if there is a finite set $S\subseteq V(G)$ so that every $S$-$\{\omega\}$ separator has order at least $k$.

Large vertex-degree at the ends together with large degree at the vertices ensures the existence of grid minors~\cite{InfExt}, and of highly connected subgraphs~\cite{hcs}, but it is not strong enough for forcing (topological) minors. This can be seen by inserting the edge set of a spanning path at each level of a large-degree tree, in a way that the obtained graph is still planar.

The reason that the vertex-degree fails to force large complete minors is  it only gives information about the sizes of vertex-separators $S_i$ `converging' to the end in question. But imagine we wish to `cut off' an end. Then the information we need is not the size  of the $S_i$, but the average amount of edges that the vertices in $S_i$ send `into' the graph. This idea is made precise in the definition of the relative degree for locally finite graphs as given in the introduction.

We defined the relative degree for locally finite graphs in the introduction as  $\inf_{(H_i)}\liminf_{i} ({|\partial_e H_i|}/{|\partial_v H_i|})$. Let us remark that instead, we might have defined the relative degree as $\inf_{(H_i)}\liminf_{i} (({|\partial_e H_i|+|E(G[\partial_v H_i])|}/{|\partial_v H_i|})$. This change, and the corresponding change for non-locally finite graphs, would alter the relative degree of the ends, but our proofs would go through in the same way. 
However, we believe that our definition is more natural.

Let us quickly evaluate a possible alternative definition, that as first sight might seem equally plausible as ours but less complicated. Consider the ratio $d_e(\omega)/d_v(\omega)$ of the edge- and the vertex-degree of some end $\omega$. If we defined this ratio as the degree of $\omega$ then large degrees at vertices and ends do not force large complete minors. For an example, see~\cite{InfExt}.

\subsection{The role of the separation property}

Let us now explain the reason for requiring the vertex-boundaries $\partial_v H_i$ of the graphs $H_i$ to be minimal  $\partial_v H_{i-1}$--$\Omega(H_i)$ separators (we shall call this property the {\em separation property}). First, let us show that a very similar condition, namely requiring the $\partial_v H_i$ to be minimal  $\partial_v H_{i-1}$--$\{ \omega \}$ separators, is too weak. Even  locally finite graphs can satisfy a degree condition thus modified and still not contain any large complete minors. In order to see this, consider the following example.

\begin{figure}
\centering
\includegraphics[width=0.7\textwidth]{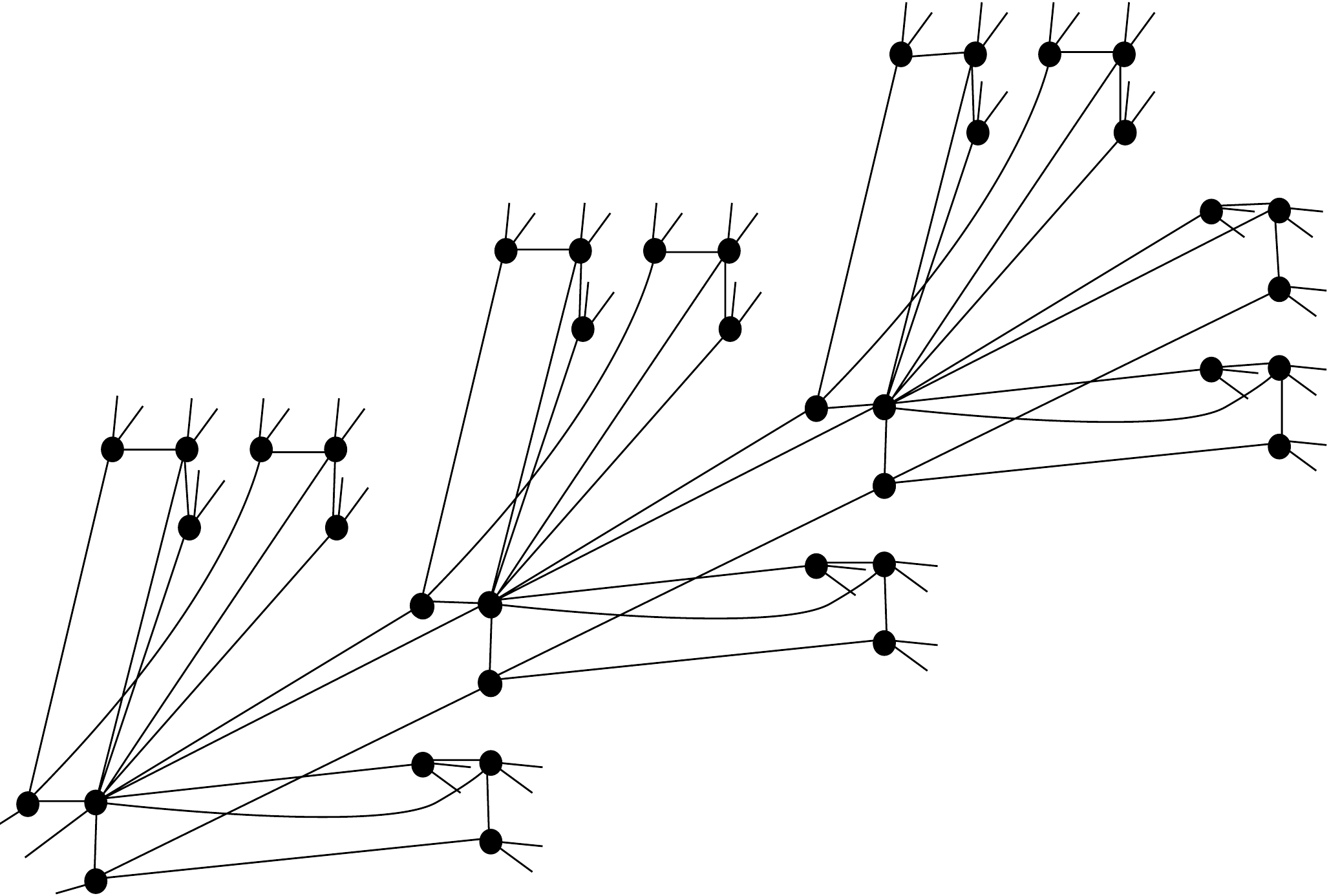} 
\caption{The graph from Example~\ref{haken}, with $r=5$.}
\end{figure}

\begin{example}\label{haken}
We start out with the $r$-regular (infinite) tree $T_r$ with levels $L_i$. For each of its vertices, we colour half of the edges going to the next level blue, the other half red. Then replace each vertex $v\in V(T_r)$ with a path $x_vy_vz_v$. For each $vw\in E(T_r)$, with $v$ from the $i$th and $w$ from the $(i+1)$th level of $T_r$, say, we do the following. If $vw$ is blue, then add the edges $x_vx_w$, $y_vy_w$ and $y_vz_w$, and if $vw$ is red, then add the edges $y_vx_w$, $y_vy_w$ and $z_vz_w$. \\
The resulting graph $\Gamma_r$ is locally finite, and it is easy to see that  $\Gamma_r$ does not contain any complete minor of order greater than $4$ using the fact that the vertices of any such minor can not be separated by any $2$-separator of $\Gamma_r$.\\
Note that each end of  $\Gamma_r$ has relative degree $1$. In fact, let $\omega\in\Omega(\Gamma_r)$, and let $\omega'\in\Omega(T_r)$ be the corresponding end of $T_r$. Suppose $v_1v_2v_3\ldots\in\omega'$ with $v_1$ being the root of $T_r$.
For $v\in V(T_r)$, let $H_v$ be the subgraph of $\Gamma_r$ that is induced by all vertices $x_w,y_w,z_w$ such that $w$ lies in the upper closure $\lceil v \rceil$ of $v$ in the tree order of $T_r$. Then $(H_{v_i})_{i\in\N}\rightarrow\omega$, and thus $\df (\omega)=1$.\\
However, if we would require the $\partial_v H_i$ to be minimal  $\partial_v H_{i-1}$--$\{ \omega \}$ separators instead of  minimal  $\partial_v H_{i-1}$--$\Omega^G(H_i)$ separators, then the ends of $\Gamma_r$ had relative degree $(r+2)/2$. In fact, as  the $\partial_v H_{v_i}$ are not minimal  $\partial_v H_{i-1}$--$\{\omega \}$ separators, the sequence $(H_{v_i})_{i\in\N}$ is no longer taken into account when we calculate the relative degree of $\omega$. It is not difficult to see that the relative degree would then be determined by the sequence $(K_{v_i})_{i\in\N}$, which we define now. If $v_iv_{i+1}$ is blue, then  $K_{v_i}$ is the subgraph of $\Gamma_r$ that is induced by $\{ x_{v_i}, y_{v_i}\}\cup\{ x_w,y_w,z_w: w\in\lceil v \rceil, v\in L_{i+1},v_iv\text{ is a blue edge}\}$. Otherwise,  let $K_{v_i}$ be induced by $\{ x_{v_i}, y_{v_i}\}\cup\{ x_w,y_w,z_w: w\in\lceil v \rceil, v\in L_{i+1},v_iv\text{ is a red edge}\}$. As ${|\partial_e K_i|}/{|\partial_v K_i|}=(r+2)/2$  for all $i>1$, the end $\omega$ would thus have relative degree $(r+2)/2$.
\end{example}

On the other hand, if we totally gave up the requirement that the $\partial_v H_i$ are minimal separators, then basically every end of every graph would have relative degree $1$. 
 To see this, write $(H_i)_{i\in\N}\rightsquigarrow\omega$ if $(H_i)_{i\in\N}$ is an infinite sequence of distinct regions of $G$ with $H_{i+1}\subseteq H_i-\partial_v H_i$ such that $\omega$ has a ray in $H_i$  for each $i\in\N$. Let  $(H_i)_{i\in\N}$ with $(H_i)_{i\in\N}\rightsquigarrow\omega$, and let $v_i\in\partial_v H_{3i}$ for $i\in\N$. Then the $v_i$ do not have common neighbours. We construct a sequence $(H'_j)_{j\in\N}$ with $H'_0:=H_0$, and, for $j>0$, we let $H'_j:=H_{i_j}-V_j$ where $i_j$ is such\footnote{For instance set $i_j:=\max\{dist(v,w)|v\in\partial_vH_0, w\in\partial_vH'_{j}\}+1$.} that $H_{i_j}\subseteq H'_j-\partial_v H'_j$, and $V_j$ consists of $j|\partial_e H_{i_j}|$ vertices $v_i$ with $i\geq i_j$. Then $(H'_j)_{j\in\N}\rightsquigarrow\omega$, and 
\[
 \liminf_{j\rightarrow\infty}\frac{|\partial_e H'_j|}{|\partial_v H'_j|}=\liminf_{j\rightarrow\infty}\frac{|\partial_e H_{i_j}|+\sum_{v\in V_j}d(v)}{|\partial_v H_{i_j}|+\sum_{v\in V_j}d(v)}=1.
\]
This shows that the additional condition that $\partial_v H_{i+1}$ is an $\subseteq$-minimal $\partial_vH_i$--$\Omega^G (H_{i+1})$ separator is indeed neccessary for our definition of the relative degree to make sense.  

For our notion, every integer, and also $\aleph_0$, appears as the relative degree of an end in some locally finite graph (larger cardinals only appear in non-locally finite graphs). Indeed, let $k\in\N$, and let $G$ be obtained from the disjoint union of $\aleph_0$ many copies $K_i$ of $K^k$ by adding all edges between $K_i$ and $K_{i+1}$, for all $i$. Suppose $(H_i)_{i\in\N}\rightarrow\omega$ for the unique end $\omega$ of $G$. Then by the separation property, we can conclude inductively that each $\partial_vH_i$ is contained in some $K_{i_j}$. Thus $\df (\omega)=k$. A similar example can be constructed to find an end of relative degree $\aleph_0$.

\subsection{From locally finite to arbitrary graphs}

In arbitrary infinite graphs, we have to face the problem that there might be vertices dominating our end $\omega$, and hence no sequence of subgraphs $H_i$ can satisfy $H_{i+1}\subseteq H_i-\partial_v H_i$. Our way out of this dilemma was to delete the dominating vertices temporarily, find the sequences as above and calculate the corresponding infimum, and then add $|Dom(\omega)|$ to the relative degree. 

One might think that alternatively, we might have weakened our requirements on the sequences of $\omega$-regions $H_i$. For instance we might be satisfied with them obeying $H_{i+1}\subseteq H_i-(\partial_v H_i-Dom(\omega))$. We then should also require that $H_i-Dom(\omega)$ is connected, as otherwise our sequence may `converge' to more than one end. And, if no sequence as just described would exist, we would set $\tilde d_{e/v}(\omega):=\infty$. This alternative definition would have the advantage that the contribution of a dominating vertex, in terms of outgoing edges, to the edge-boundary of an $\omega$-region is counted. 

However, the approach does not allow for an extension of Theorem~\ref{thm:minor}. The problem are vertices that dominate more than one end. Consider an infinite $r$-regular tree to which we add one vertex that is adjacent to all other vertices. The ends of this graph have infinite degree in the sense just discussed (and relative degree $2$), but of course, no $K^k$-minor for large $k$, no matter how large $r$ is with respect to $k$. We can give a similar example for a graph with only two ends.\footnote{Take two copies of the $r$-regular tree add a spanning path in each level of each of the two trees and a new vertex adjacent to all other vertices.}

\section{\hskip-.18cm Dominating vertices and topological $K^{k}$-minors}

This section provides some results about dominating vertices that will be needed later on.
First, we give a useful characterization of dominating vertices.

\begin{lemma}\label{lem:domi}
In any graph, a vertex $v$ dominates an end $\omega$ if and only if there is no finite  $v$--$\omega$ separator.
\end{lemma}
\begin{proof}
For the forward direction, note that every finite set of vertices intersects only a finite number of the infinitely many $v$--$V(R)$ paths, where $R\in \omega$ is dominated by $v$. Hence $v$ and $\omega$ cannot be finitely separated. 

For the backward direction we inductively construct a set of   $v$--$V(R)$ paths, where $R$ is any ray in $\omega$. At each step we use the fact that all paths constructed so far form a finite set which (without $v$) does not separate $v$ from $\omega$. Hence we can always add a new $v$--$V(R)$ path, which disjoint from all the others (except in $v$). 
\end{proof}

We now show that for an end $\omega$ which is dominated by only finitely many vertices, the graph $G_\omega$ contains sequences of subgraphs converging to $\hat\omega$. 

\begin{lemma}\label{lem:sequ}
Let  $G$ be a graph, and  let $\omega\in\Omega(G)$ with $|Dom (\omega)|<\aleph_0$. Then $G_\omega$ contains a sequence $(\hat H_i)_{i\in\mathbb N}\rightarrow\hat\omega$. 
\end{lemma}

\begin{proof}
 We define a sequence of disjoint finite sets $S_i\subseteq V(G_\omega)$, starting with any finite non-empty set $S_1$. For $i>1$ and for each $v\in S_{i-1}$ let $S_v^{i-1}$ be a finite $v$--$\hat\omega$ separator in $G_\omega$ (note that such a separator exists by Lemma~\ref{lem:domi}). Set $$\tilde S_i:=\bigcup_{v\in S_{i-1}}S_v^{i-1}\setminus S_{i-1}.$$ Then $\tilde S_i$ separates $S_{i-1}$ from $\hat\omega$. In fact, otherwise there would be a ray $R\in\hat\omega$ with only its first vertex $v$ in $S_{i-1}$, and disjoint from $\tilde S_i$. But $R$ must meet $S_v^{i-1}$, a contradiction.
 
 Choose $S_i\subseteq \tilde S_i$ minimal such that it separates $S_{i-1}$ from $\hat\omega$. Now for $i\in\mathbb N$, let $K_i$ be the  component of $G_\omega-S_i$ that contains a ray of $\hat\omega$ (since $S_i\cup Dom(\omega)$ is finite there is a unique such component $K_i$). Let $\hat H_i$ be the subgraph of $G_\omega$ that is induced by $S_i$ and $K_i$. Then, by the choice of $S_i$, we find that $S_i$ is a minimal $S_{i-1}$--$\Omega^G(\hat H_i)$ separator. Hence, $(\hat H_i)_{i\in\N}\to\hat\omega$, as desired.
\end{proof}

Next, we will see that our desired minor is easy to find whenever there are enough vertices dominating the same end.

\begin{lemma}\label{lem:minor}
Let $k\in\mathbb N$, let $G$ be a graph and let $\omega\in\Omega(G)$. If $|Dom(\omega)|\geq k$, then $K^{k}\preceq_{top} G$. 
\end{lemma}

Lemma~\ref{lem:minor} follows at once from Lemma~\ref{lem:minor2} below. The {\em branching vertices} of a subdivision are those vertices that did not arise from subdividing edges.

\begin{lemma}\label{lem:minor2}
Let $k\in\mathbb N$, let $G$ be a graph,  let $\omega\in\Omega(G)$, and let $S\subseteq Dom(\omega)$ with $|S|= k$. Then $G$ contains a subdivision $TK^{k}$ of $K^{k}$ whose branching vertices are in $S$. 
\end{lemma}

\begin{proof}
We use induction on $k$, the base case $k=0$ is trivial. So suppose $k\geq 1$. Then, let $S\subseteq Dom(\omega)$ be a set of size $k$, and let $s\in S$. By the induction hypothesis, $G-s$ contains a subdivision $TK^{k -1}$ of $K^{k -1}$ with branching vertices in $S':=S\setminus\{s\}$. 

Successively we define sets $\mathcal P_i$ of $s$--$S'$ paths in $G$ which are disjoint except in $s$. We start with $\mathcal P_0:=\emptyset$. For $i>0$, suppose there is a vertex $v\in S'$ which is not the endpoint of a path in $\mathcal P_{i-1}$. Then, $S_i:=S\cup \bigcup_{P\in\mathcal P_{i-1}}V(P)$ is finite, and $G-S_i$ has a unique component $C_i$ which contains rays of $\omega$.  Since Lemma~\ref{lem:domi} implies that neither $s$ nor $v$ can be separated from $\omega$ by a finite set of vertices, both $s$ and $v$ have neighbours in $C_i$. Hence there is an $s$--$v$ path $P_i$ that is internally disjoint from  $S_{i}$. Set $\mathcal P_i:=\mathcal P_{i-1}\cup \{P_i\}$. The procedure stops after step $k-1$, when all vertices of $S'$ are connected to $s$ by a path in $\mathcal P_i$. This gives the desired subdivision $TK^{k}$.
\end{proof}

With a very similar proof,\footnote{We construct the $TK^{\aleph_0}$ step by step, adding one branching vertex plus the corresponding paths at a time.  In each step, the finiteness of the already constructed part ensures the existence of an unused dominating vertex $s$ of $\omega$ and enough paths to connect $s$ to the already defined branching vertices.} we also get the following statement (which will not be needed in what follows):

\begin{lemma}\label{lem:minor3}
Let $G$ be a graph and let $\omega\in\Omega(G)$. If $|Dom(\omega)|\geq \aleph_0$, then $K^{\aleph_0}\preceq_{top} G$. 
\end{lemma}

We finish this section with one more basic lemma. This lemma implies that removing a finite part of a graph, or even an infinite part with a finite vertex-boundary, will not alter the relative degree of the remaining ends. 

\begin{lemma}\label{lem:grad}
Let $G$ be a graph, let $\omega \in \Omega (G)$, and let $G'$ be an induced subgraph of $G$ such that $\partial_v^G G'$ is finite, and such that $G'$ has an end $\omega'$ that contains rays of $\omega$. Then $\df (\omega')=\df (\omega)$.
\end{lemma}

\begin{proof}
First of all, observe that since $\partial_v^G G'$ is finite, it follows that $Dom(\omega)\subseteq V(G')$. Hence, we only need to show that $\inf_{(\hat H_i)_{i\in\N}}\liminf_{i}({|\partial_e \hat H_i|}/{|\partial_v \hat H_i|})$, is the same for sequences $(\hat H_i)_{i\in\N}$ in $G_\omega$ and for sequences $(\hat H_i)_{i\in\N}$  in $G'_{\omega'}$.

For this, note that every sequence $(\hat H_i)_{i\in\N}\rightarrow \hat\omega$ in $G_\omega$ has a subsequence $(\hat H_i)_{i\geq i_0}\rightarrow \hat\omega'$ in $G'_{\omega'}$ (it suffices to take $i_0:=\max\{dist(v,w):v\in \partial_v \hat H_0, w\in \partial_v G'\}+1$). Moreover, for every sequence $(
\hat H'_j)_{j\in\N}\rightarrow \hat\omega'$ in $G'_{\omega'}$, there is an index $ j_0$ such that $\partial_v \hat H'_j\cup N(\hat H'_j)$ is the same set in $G_{\omega}$ and in $G'_{\omega}$ (for instance, take $j_0:=\max\{dist(v,w):v\in \partial_v\hat H'_0, w\in \partial_v G'\}+1$). Hence, the relative degrees of $\omega$ and $\omega'$ are the same.
\end{proof}

\section{Proof of Theorem~\ref{thm:cntble}}

In this section we prove our main result, Theorem~\ref{thm:cntble}.
We start by showing how to find, for a fixed end $\omega$ of some graph $G$, an $\hat\omega$-region $\hat H$ of $G_\omega$ that has an acceptable average degree into $G-\hat H$.

\begin{lemma}\label{cutoffnew1}
  Let $G$ be a graph, let $\omega\in\Omega (G)$ with $\df (\omega)> m$ for some $m\in\mathbb Q$, and let $S\subseteq V(G)$ be finite. If $|Dom(\omega)|<\aleph_0$, then $G_\omega$ has a $\hat\omega$-region $\hat H$ such that 
  \begin{enumerate}[(a)]
\item $S\cap V(\hat H) =\emptyset$, and
\item  $\frac{|\partial_e \hat H|}{|\partial_v \hat H|}> m-|Dom(\omega)|$. 
  \end{enumerate}
 \end{lemma}

\begin{proof}
By Lemma~\ref{lem:sequ}, there is a sequence $(\hat H_i)_{i\in\mathbb N}\rightarrow\hat\omega$ in $G_\omega$. 
Let $i_0=\max\{dist_{G_\omega}(v,w):v\in S\setminus Dom (\omega), w\in\partial_v \hat H_0\}$+1. Then $S\cap V(\hat H_i)=\emptyset$ for all $i\geq i_0$.
As $\df (\omega)> m$, there is a $j_0 \geq i_0$ such that $({|\partial_e \hat H_j|}/{|\partial_v \hat H_j|})> m- |Dom(\omega)|$ for all $j\geq j_0$. 
\end{proof}

We now apply Lemma~\ref{cutoffnew1} repeatedly to  the ends of any suitable fixed countable subset of $\Omega(G)$. Lemma~\ref{lem:grad} will ensure that the relative degree of the ends is not disturbed by what has been cut off earlier.

\begin{figure}[h]
\centering
\scalebox{0.4}{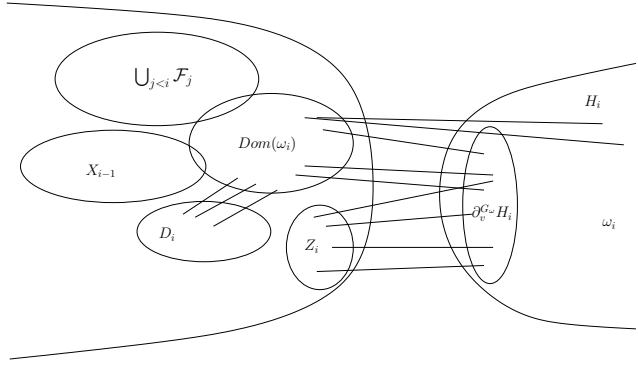}
\caption{Construction of the graph $G_i$ in Lemma \ref{lem:Gi}.}
\end{figure}

\begin{lemma}\label{lem:Gi}
Let $m\in \mathbb Q$, let $G$ be a graph,  let  $X\subseteq V(G)$ be finite and for $i\in\N$, let $\omega_i\in\Omega(G)$ with $|Dom(\omega_i)|<\aleph_0$. If $\delta^{V,\Omega}(G)>m$,  then $G$ contains induced subgraphs $G_i$, finite sets $X\subseteq X_i\subseteq V(G_i)$, and finite sets $\mathcal F_i$ of pairwise disjoint finite subsets of $V(G_i)$ such that for each $i\in\N$
\begin{enumerate}[(A)] 
\item $G_{i+1}\subseteq G_i$, $X_i\subseteq X_{i+1}$,  and $\mathcal F_{i+1}\supseteq \mathcal F_i$,\label{A}
\item $X_i\cup \bigcup\mathcal F_i\subseteq V(G_i)$, \label{u}
\item there is a family $\mathcal H_i=\{ H_F:F\in\mathcal F_i\}$ of disjoint connected subgraphs such that $V(H_F)\cap V(G_i)= F$ for each $F\in\mathcal F_i$ and such that $G=G_i\cup\bigcup\mathcal H_i$,\label{H_F}
\item $\partial_v^G G_i$ is finite,\label{rand}
\item the average degree of $F$ into $F\cup X_i$ is  $>m-|Dom(\omega_i)|$, for each $F\in \mathcal F_i$, \label{avdeg}
\item $d_{X_i}(v)>m$ for all $v\in \partial_v^G G_i\setminus \bigcup\mathcal F_i$, and\label{mindeg}
\item $\omega_i$ has no rays in $G_i$.\label{norays}
\end{enumerate}
\end{lemma}

\begin{proof}
Set $G_0:=G$,  $X_0:=X$, and  $\mathcal F_0:=\emptyset$. Now, for
 $i\geq 1$ we do the following. If $\omega_i$ has no ray in $G_{i-1}$, then we set $G_i:=G_{i-1}$, $X_i:=X_{i-1}$ and $\mathcal F_i:=\mathcal F_{i-1}$, which ensures all the desired properties (as they hold for $i-1$). 

So suppose $\omega_i$ does have a ray  in $G_{i-1}$. Then let $D_i$ be  a  finite set of vertices of $G_{i-1}$ so that each dominating vertex of $\omega_i$ in $G_{i-1}$ has degree greater than $m$ into $D_i$. (This is possible since by assumption there only finitely many vertices dominating $\omega_i$.)

Observe that because of~\eqref{rand} we may apply Lemma~\ref{lem:grad} to obtain that all ends of $G_{i-1}$ have relative degree $>m$. Hence Lemma \ref{cutoffnew1} applied to $G_{i-1}$ and the finite set 
$$S:=X_{i-1}\cup D_i\cup \bigcup\mathcal F_{i-1}$$
yields an $\hat\omega_i$-region $\hat H_i$ of $G_\omega$.
We set $\mathcal F_i:=\mathcal F_{i-1}\cup \partial_v \hat H_i $ and $G_i:=G_{i-1}-(\hat H_i -\partial_v \hat H_i)$. Choose  a finite subset $Z_i$ of $V(G_{i})$ such that $ \partial_v \hat H_i $ has average degree $>m-|Dom(\omega_i)|$ into $Z_i\cup \partial_v \hat H_i $ and set $X_i:=X_{i-1}\cup D_i\cup Z_i$.

Then, conditions~\eqref{A},~\eqref{rand} and~\eqref{norays} are clearly satisfied  for step $i$, as they hold for step $i-1$. Conditions~\eqref{u} and~\eqref{H_F} for $i$ follow from Lemma~\ref{cutoffnew1}(a), and from~\eqref{u} and~\eqref{H_F} for $i-1$.
Condition~\eqref{avdeg} follows from Lemma~\ref{cutoffnew1}(b)  and~\eqref{avdeg} for $i-1$.

 Finally, for~\eqref{mindeg} suppose that $v\in \partial_v^G G_i\setminus \bigcup\mathcal F_i$. If $v\in \partial_v^G G_{i-1}$ then~\eqref{mindeg} for $i$ follows from~\eqref{mindeg} for $i-1$. Otherwise, $v$ dominates  $\omega_i$. Then by construction, $v$ has sufficiently many neighbours in $D_i\subseteq X_i$. 
 \end{proof}

If $G$ has only countably many ends, then the procedure just described can be used to cut off all ends:

\begin{lemma}\label{lem:cutoffallends}
Let $k\in\N$, let $m\in \mathbb Q$, let $G$ be a graph with $|\Omega(G)|\leq\aleph_0$, and let  $X\subseteq V(G)$ be finite. Suppose $|Dom(\omega)|<k$ for all ends $\omega\in\Omega (G)$. If $\delta^{V,\Omega}(G)>m$,  then $G$ has an induced subgraph $G'$ and a set $\mathcal F$ of finite pairwise disjoint vertex sets such that
\begin{enumerate}[(i)]
\item $X\cup \bigcup\mathcal F\subseteq V(G')$, \label{uu}
\item there is a family  $\{ H_F:F\in\mathcal F\}$ of disjoint connected subgraphs of $G$ such that $V(H_F)\cap V(G')= F$ for each $F\in\mathcal F$, \label{H_FF}
\item for each vertex $v\in V(G')$ of degree  $\leq m$ there is an $F\in\mathcal F$ so that $v\in F$ and the average degree of $F$ in $G'$ is $>m-k+1$, and\label{deg}
\item every ray of $G$ has only finitely many vertices in $V(G')$.\label{rays}
\end{enumerate}
\end{lemma}

\begin{proof}
Let $\omega_1,\omega_2,\omega_3,\ldots$ be a (possibly repetitive) enumeration of $\Omega(G)$.
Apply Lemma~\ref{lem:Gi}, and then set $G':=\bigcap_{i\in\N} G_i$ and $\mathcal F:=\bigcup_{i\in\N}\mathcal F_i$. We claim that $G'$ and $\mathcal F$ are as desired. Indeed, properties~\eqref{A} and~\eqref{u} imply property~\eqref{uu}, and property~\eqref{H_F} together with~\eqref{A} implies property~\eqref{H_FF}. 

For property~\eqref{deg} observe that~\eqref{A} and~\eqref{u} imply that $X_i\subseteq V(G')$ for all $i\in\mathbb N$. Now~\eqref{deg} follows from~\eqref{avdeg} and~\eqref{mindeg} together with the assumption that $\delta^{V,\Omega}(G)>m$.

In order to see~\eqref{rays}, suppose that $R$ is a ray of $G$ that has infinitely many vertices in $G'$. Say $R\in  \omega_j$. Then by~\eqref{A} for $j$, the ray $R$ has infinitely many vertices in $G_j$. So, as $\partial_v^GG_j$ is finite by~\eqref{rand}, $R$ has a subray in $G_j$, a contradiction to~\eqref{norays} for~$j$.
\end{proof}

We are now almost ready to prove our main theorem. We will make use of a standard tool from infinite graph theory, K\H onig's infinity lemma.

\begin{lemma}\label{koenig}$\!\!${\bf\cite{DBook}}
Let $V_1, V_2, V_3,\ldots$ be disjoint finite non-empty sets, and let  $G$ be a graph on their union. Suppose that for all $i\in\N$, each  vertex of $V_{i+1}$ has a neighbour in $V_i$. Then $G$ has a ray $v_1v_2v_3\ldots$, with $v_i\in V_i$, for each $i\in\N$.
\end{lemma}

Let us now prove Theorem~\ref{thm:cntble}.

\begin{proof}[Proof of Theorem~\ref{thm:cntble}]
Suppose that  $K^{k}$ is not a topological minor of $G$. Then by Lemma~\ref{lem:minor}, $|Dom(\omega)|<k$ for all $\omega\in\Omega (G)$. Let $u\in V(G)$ and let $G'$ be the subgraph we obtain from Lemma~\ref{lem:cutoffallends} applied to $G$ and $X:=\{ u\}$, and let $\mathcal F$ be the corresponding set of disjoint finite vertex sets. 

For $i\in\mathbb N$, we shall successively define finite sets $S_i$, with $S_i\subseteq S_{i+1}$. We start with setting $S_0:=\emptyset$ and $S_1:=\{u\}$ if $u\notin\bigcup \mathcal F$, or $S_1:=F_u$ if there is an $F_u\in\mathcal F$ with $u\in F_u$ (by the disjointness of the sets in $\mathcal F$, there is at most one such $F_u$). Note that $S_0\subseteq S_1\subseteq V(G')$ because of Lemma~\ref{lem:cutoffallends}~\eqref{uu}. 

Our sets $S_i$ will have the following properties for $i\geq 1$:
\begin{enumerate}[(I)]
\item the average degree of the vertices of $S_{i-1}$ in $G'[S_i]$ is $>m-k+1$, \label{a}
\item the average degree of the vertices of $S_i\setminus S_{i-1}$ in $G'$ is $>m-k+1$, and\label{b}
\item for each $F\in\mathcal F$ with $F\cap S_i\neq \emptyset$ there is a $j\leq i$ such that $F\subseteq S_j\setminus S_{j-1}$.\label{c}
\end{enumerate}

For $i=1$, property~\eqref{a} holds trivially, and~\eqref{b} is satisfied because of Lemma~\ref{lem:cutoffallends}~\eqref{deg}. By the choice of $S_1$ and since the $F\in\mathcal F$ are disjoint, also~\eqref{c} holds.

Now, for $i\geq 2$ we choose a finite subset $X_i$ of the neighbourhood of $S_{i-1}\setminus S_{i-2}$ in $G'-S_{i-1}$ so that  the average degree of the vertices in $S_{i-1}$ in the graph $G'[S_{i-1}\cup X_i]$ is at least $m-k+1$. Such a choice is possible by~\eqref{a} and~\eqref{b} for $i-1$.

Let $Y_i$ denote the union of all $F\in\mathcal F$ that contain some $v\in X_i$ . Note that  $Y_i$ is finite since the $F$ are all disjoint and because $X_i$ is finite.  Then set $S_i:=S_{i-1}\cup X_i \cup Y_i$. Our choice of the $S_i$ clearly satisfies conditions~\eqref{a}, \eqref{b}  and~\eqref{c}. This finishes our definition of the sets $S_i$.

First suppose  that $S_i\neq S_{i-1}$ for all $i\in\N$. Then, for each $i\in\N$, let $V_i$ be  obtained from $S_i\setminus S_{i-1}$ by collapsing each $F\in\mathcal F$ with $F\subseteq S_i\setminus S_{i-1}$ to one vertex $v_F$ (which will be adjacent to all neighbours of $F$ outside $F$). So, for all $i\in\N$, each vertex of $V_{i+1}$ has a neighbour in $V_i$, and therefore, we may apply K\H onig's infinity lemma (Lemma~\ref{koenig}) to the sets $V_i$ in order to find a ray $R$ in $\bigcup_{i\in\N} V_i$. We use \eqref{c} and Lemma~\ref{lem:cutoffallends}~\eqref{H_FF}  to expand $R$ to a ray $R'$ in $G$. As $R'$ has infinitely many vertices in $G'$, this establishes a contradiction to Lemma~\ref{lem:cutoffallends}~\eqref{rays}.

So we may assume that there is an $i\in \N$ such that $S_i=S_{i-1}$. Then,  by~\eqref{a}, $G'[S_i]$ is a finite graph of average degree $>m-k+1$, which is as desired. 
 \end{proof}

 \bibliographystyle{plain}
 \bibliography{graphs}

\small

\end{document}